\documentclass{article}
\usepackage{amsmath,graphicx,mlspconf}

\usepackage{amsmath,cite,xspace}
\usepackage{graphicx}
\usepackage{epsfig}
\usepackage{epstopdf}
\usepackage{bm}
\usepackage{color}
\usepackage{amssymb}
\usepackage{hyperref}


\usepackage{mathtools}

\newcommand{\bp}{{\bf p}}

\newcommand{\bv}{{\bf v}}

\newcommand{\cV}{{\cal V}}

\graphicspath{{./Figures/}} \DeclareGraphicsExtensions{.eps,.pdf,.jpg,.mps,.png}
\begin{document}

\title{Learning Warm-Start Points for AC Optimal Power Flow } 

\author{Kyri Baker, 
\thanks{K. Baker is with the College of Engineering and Applied Science at the University of Colorado, Boulder. Email: kyri.baker@colorado.edu.}}
\name{Kyri Baker}
\address{University of Colorado Boulder}

\maketitle

\begin{abstract}
A large amount of data has been generated by grid operators solving AC optimal power flow (ACOPF) throughout the years, and we explore how leveraging this data can be used to help solve future ACOPF problems. We use this data to train a Random Forest to predict solutions of future ACOPF problems. To preserve correlations and relationships between predicted variables, we utilize a multi-target approach to learn approximate voltage and generation solutions to ACOPF problems directly by only using network loads, without the knowledge of other network parameters or the system topology. We explore the benefits of using the learned solution as an intelligent warm start point for solving the ACOPF, and the proposed framework is evaluated numerically using multiple IEEE test networks. The benefit of using learned ACOPF solutions is shown to be solver and network dependent, but shows promise for quickly finding approximate solutions to the ACOPF problem.

\end{abstract}

\section{Introduction and Motivation} \label{sec:intro}
Grid operators repeatedly solve optimal power flow across transmission networks, multiple times throughout the day, every day of the year, for decades, to ensure that the grid is operating reliably and safely. They do not stop for holidays, weekends, or even birthdays, adhering to the reliability standards defined by the Energy Policy Act of 2005 and enforced by the North American Electric Reliability Corporation (NERC). The utility of this data has been recognized by power system operators in the past, and has been used since the late 1960's to solve power system state estimation problems \cite{Schweppe70, Wu90}, whose goal is to estimate the complex nodal voltages within a power system. Inspired partially by these data-driven approaches, a natural extension of this concept is to now leverage data for estimating the solution of AC \emph{optimal} power flow problems.

High-voltage transmission grids are typically modeled as a set of $N$ buses $\cal{N}$, generators $\cal{G} \subseteq \cal{N}$ where $|\cal{G}| >$ 0, and a bus admittance matrix describing the topology of the network, $Y = G + jB \in \mathbb{C}$. The AC power flow (ACPF) equations, which solve for complex voltages $|v_m|e^{j\theta_{v_m}} \in \mathbb{C}$ throughout the network given injected complex powers $s_m = p_m + jq_m \in \mathbb{C}$, represent the AC steady-state conditions of the grid. These prototypical AC power flow equations in rectangular coordinates are given in terms of the voltage magnitude at each bus $m \in \cal{N}$ as $|v_m|$, voltage angle $\theta_{v_m}$, net real power $p_m$, and net reactive power $q_m$ as

\begin{subequations} \label{eqn:ACPF}
\begin{align} 
p_m = |v_m| \sum_{l \in \cal{N}}v_l(G_{ml}\cos(\theta_{ml}
) + B_{ml}\sin(\theta_{ml})),\label{eqn:PF_constr}\\
q_m = |v_m| \sum_{l \in \cal{N}}v_l(G_{ml}\sin(\theta_{ml}
) - B_{ml}\cos(\theta_{ml})),\label{eqn:PF_constr2}
 \end{align}
\end{subequations}

\noindent where the angle difference $\theta_{ml} := \theta_m - \theta_l$. Using the ACPF equations as physical constraints, the AC optimal power flow (ACOPF) problem is thus created. A solution to this problem seeks to satisfy the physical power flows while minimizing generation cost and adhering to system constraints. This problem can generally be written as
\begin{subequations} \label{eqn:ACOPF}
\begin{align} 
 \, \,  &\min_{\substack{\bv \in \cV, \bp_g}}\hspace{.2cm} \sum_{j \in \cal{G}} a_j p_{g,j}^2 + b_j p_{g,j} + c_j \\
&\mathrm{subject\,to:~} \, \eqref{eqn:PF_constr}, \eqref{eqn:PF_constr2}\\
&\hspace{2cm} \, \underline{p}_{g,j} \leq p_{g,j} \leq \overline{p}_{g,j}, \forall j \in \cal{G} \\
&\hspace{2cm} \, \underline{q}_{g,j} \leq q_{g,j} \leq \overline{q}_{g,j}, \forall j \in \cal{G} 
 \end{align}
\end{subequations}

\noindent where coefficients $a_j$, $b_j$, and $c_j$ represent the cost of generator $j$, $\cal{V}$ is the set of permissible complex voltages at each node in $\cal{N}$, $p_{g,j}$ ($q_{g,j}$) is the controllable active (reactive) power output of generator $j$, and $\underline{p}_{g,j}$ ($\underline{q}_{g,j}$) and $\overline{p}_{g,j}$ ($\overline{q}_{g,j}$) are lower and upper limits on active (reactive) power generation, respectively.

In this paper, we investigate potential new ways of solving the ACOPF problem by leveraging machine learning techniques to both expedite finding a solution of \eqref{eqn:ACOPF} and to find approximate solutions to \eqref{eqn:ACOPF} that perform better than current methods that grid operators rely on to solve this difficult problem. It is also worth noting that the method developed in this paper is not specific to ACOPF problems, and learning approaches can be explored for obtaining better initial guesses in many applications of optimization.

\subsection{Challenges with solving ACOPF}
Due to the nonconvex, nonlinear nature of the ACOPF problem and constraints, finding a solution with the use of iterative numerical methods may result in a failure of the method to converge. Using warm start solutions to initialize the iterative procedure (in traditional operation, this is usually a Newton-Raphson based iterative method) of solving the power flow equations can help ensure convergence \cite{dong12}, and is typically  performed by grid operators, generally by using the previous timestep's ACOPF solution to intialize the current timestep, using a flat start, or by first solving an easier, linearized version of the ACOPF equation under the presumption that these solutions will be relatively close. However, a flat start may not ensure a solution to \eqref{eqn:ACPF}; first solving a convex approximation of \eqref{eqn:ACOPF} to obtain a better initial point requires extra computational time; and using the previous timestep's solution may not be a good approximation of the current solution, especially as more intermittent renewable energy sources and inverter-interfaced assets become grid connected.


In this paper, we aim to answer two research questions that address the aforementioned challenges: \emph{How close can learning-based approaches get us to the solution of the ACOPF problem? Are purely data-driven methods worth pursuing further?}; and \emph{Can we use learning to find better warm-start initial guesses when solving the ACOPF problem?} In addition to these, using learning approaches to find successful warm start points for iterative optimization algorithms is scarcely mentioned in the literature, for power systems problems or otherwise. We hope to help further the work in this area outside of ACOPF problems.

\subsection{Previous related efforts}
The importance of finding robust solutions to \eqref{eqn:ACOPF} or finding approximate solutions quickly has been an important topic in power systems engineering since the development of the original AC power flow equations. Perhaps the most well-known method of addressing computational and convergence issues is by reformulating the ACOPF as a linearized DC optimal power flow problem (DCOPF), which involves various physical assumptions such as lossless lines, voltage magnitudes equalling $1.0$, and a lack of reactive power control or compensation \cite{Overbye04}. Other various linearizations of \eqref{eqn:ACPF} have been developed over the years for power system optimization, including decoupled and fast decoupled loadflow \cite{decoup_loadflow}, and data-driven techniques that can estimate linearizations of \eqref{eqn:ACPF} \cite{Li_linear_18, Misra_linear_18}. 

Quality and accuracy of the linearizations is highly dependent on the scenario, usually rely on model-based solutions or knowledge of the network topology, and can result in large errors between the linearized model output and true ACOPF solution. To maintain guaranteed convergence while preserving the use of the original ACPF equations, some works address the limitations of poor starting points in Newton-based methods through the use of continuation-based methods \cite{Continuation91, Sanja12, Milano09} which alleviate many issues with poorly chosen initial conditions. However, these methods are generally slow, greatly favoring one side of the robustness versus speed tradeoff. In \cite{Venzke19}, the benefits of warm-starting ACOPF problems with the solution from convex relaxations is shown from both computational and optimality perspectives. In this work, they indicated that improving the warm start point could have convergence and computational speed benefits, and that the benefits of a warm-start were solver dependent, which we also conclude in this work. 

With these limitations and challenges in the current state-of-the-art, we wish to pursue approximate solutions to the original ACOPF problem while minimizing convergence issues and pursuing a reasonable computation time. We attempt to further the work in this area by using machine learning techniques. The power of using learning for power systems engineering has become evident in recent years \cite{Zamzam19, Misra19, pred_OPF_16}. In \cite{Misra19}, a learning approach is applied to learning the active sets of the linearized DC optimal power flow problem, with promising results when the number of active sets is small. Perhaps the closest previous work to this paper is \cite{pred_OPF_16}, which looked at the potential of predicting solutions to the ACOPF problem on a small (30 bus) test network to improve computation time. The main challenge of learning OPF solutions directly is elucidated in this unpublished work: approximately 60\% of the time, constraints were violated in the optimization problem. This is also a large risk encountered in learning which constraints are binding in linearized stochastic ACOPF problems \cite{BakerTSG19}. Thus, we aim to continue and improve upon the limited work in this area by learning \emph{warm start} solutions to ACOPF, preserving the use of the original ACPF equations but improving convergence properties and computational speed, and by further analyzing the capabilities of learning to contribute to the field of optimal power flow. Towards these goals, we build an ensemble of learners using a multi-target Random Forest algorithm \cite{StatLearn}, which is a supervised machine learning method that combines the output of multiple agents (individual Decision Trees) to produce a more accurate prediction of the ACOPF solution.

\section{Learning Framework} \label{sec:framework}
In this framework, we assume no knowledge about the network topology (i.e., we do not have access to the admittance matrix or line parameters) other than the power demands, knowledge of how many generators are in the network, and how many buses are in the network. The inputs (features) to the Random Forest are simply the loads at each bus. The outputs of the model are the optimal generation values $p^*_{g,j}$ and the voltage magnitudes $|v^*|$. We also assume that we have access to thousands of previous optimal power flow solutions that we can use to train a Random Forest.

\subsection{Multi-target Random Forests} Generally, machine learning problems, be it either regression or classification, focus on situations where there are multiple independent variables (the predictors) and a single dependent variable (what we want to predict or classify). In this problem, the voltages and optimal generation values we want to predict can be correlated, and we wish to capitalize on this fact to improve the accuracy of our learning approach. This approach differs from \cite{pred_OPF_16}; in order to preserve correlations between output variables (voltages and powers), we leverage recent advancements in \emph{multi-output}, or \emph{multi-target} regression \cite{multi_output_regression}, which allows us to simultaneously predict multiple outputs (in this case, the optimal voltage magnitudes and active power generation in the ACOPF problem). 

One benefit of using ensemble methods such as Random Forests is that the individual Decision Tree regressions can be done in a distributed and completely decentralized way, unlike centralized learning approaches such as neural networks. Each Decision Tree agent performs their individual regression procedure separately, and then agents coordinate at the very end for a final quick averaging procedure across soluions, allowing for scalability towards large-scale transmission networks. In addition, the use of a Random Forest helps alleviate risks of overfitting associated with individual Decision Trees. 

\begin{figure}[t!]
\centering
    \includegraphics[width=0.45\textwidth]{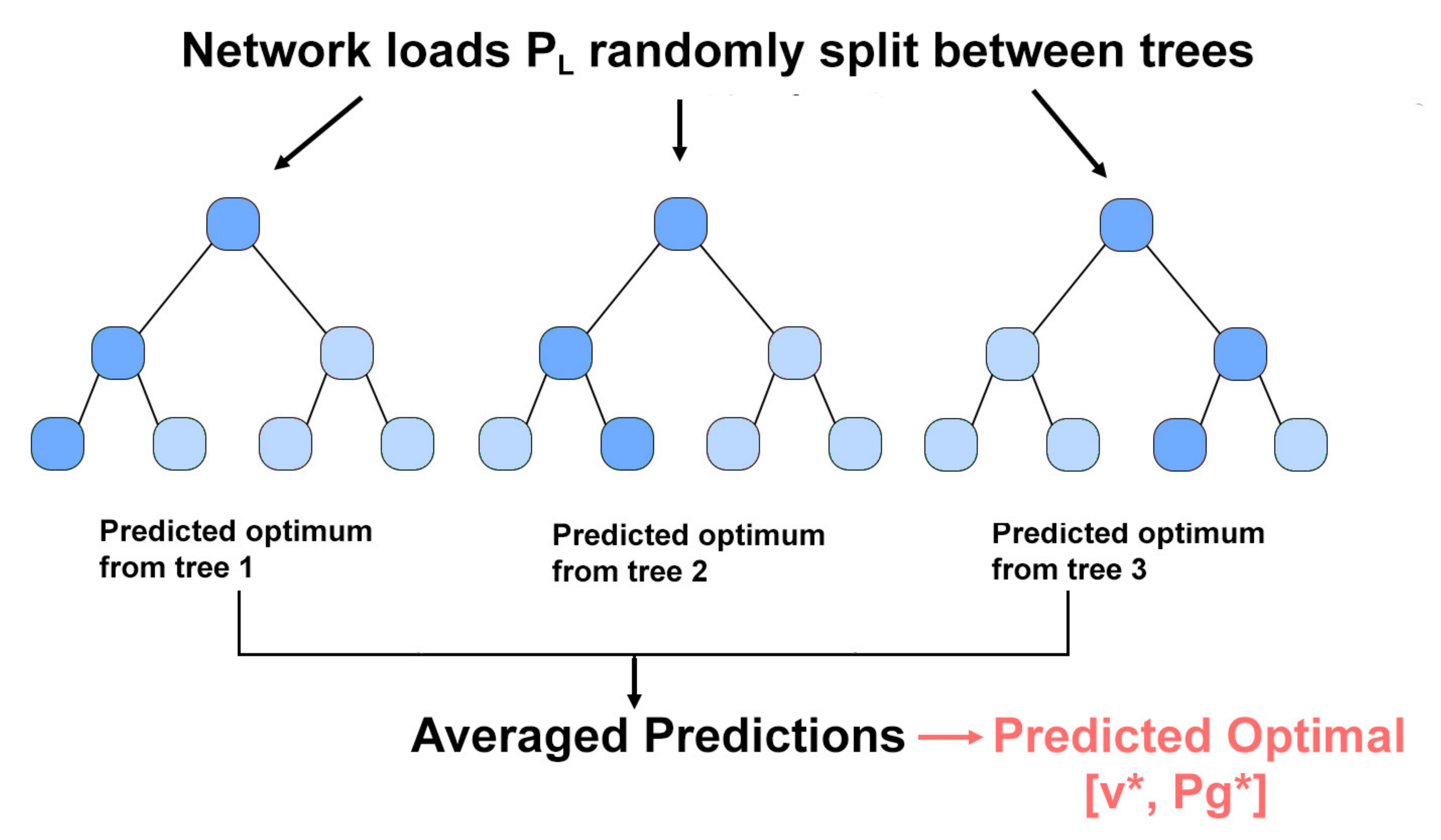}
    \caption{The multi-target Random Forest combines the outputs of multiple Decision Tree regressions (which are computed in a distributed fashion) to produce a single predicted optimal vector of voltages and generation.}
    \label{fig:RF}
\end{figure}

\subsection{Warm starts - preventing forest fires}
After the Random Forest has been trained using the ACOPF solutions in the training dataset, the model can now be used to predict optimal voltage and generation values given the network loading as an input. Naturally, even with a very low error rate in the predicted variables, the Forest does not necessarily ensure that the constraints in the original OPF problem are satisfied. This disconnect between regression and constraint satisfaction is the motivation behind our focus on learning solutions that result in better warm starts instead of learning the optimal solution directly.

Note that here we make no claims about global optimality of the ACOPF solution. Methods pursuing the global optimum of ACOPF do exist \cite{Lavaei_zero, global1} but may not be practical when solving ACOPF on a sub-second level. Generating these solutions offline and using them to train the Random Forest is an interesting approach that deserves consideration in future work, and could also be used to avoid the potential problem of having multiple training samples corresponding to a single set of network loads. Warm starts have also been shown to help achieve the global optimum in some cases \cite{Venzke19}.

\section{Numerical Experiments} \label{sec:simulations}
\subsection{System setup}
We used MATPOWER \cite{MATPOWER} to generate the training and testing datasets, and to perform ACOPF simulations with different initial starting points. The chosen IEEE test cases were solved repeatedly for uniformly distributed random perturbations of the load at each bus between $120\%$ and $80\%$ of the default load values. Solutions that resulted in infeasibilities (not enough generation capacity, etc.) were not included in the dataset. The dataset was split 80/20 into training and testing data, and passed to the \texttt{RandomForestRegressor} from Python's Scikit-learn package. The four IEEE networks considered in these experiments were the IEEE 14-bus, 57-bus, 118-bus, and 300-bus networks.\\


\noindent Recall that the number of dependent variables that we are simultaneously predicting is equal to the number of generators plus the number of buses; for the above test networks (14, 57, 118, and 300 buses), this equates to 19, 64, 172, and 369 predicted output variables, respectively. 2000 samples (runs of the ACOPF) were generated for each network (1600 training, 400 testing). MATPOWER's \texttt{MIPS} (Matpower Interior Point Solver), based on a primal-dual interior point method, was used for the iteration comparison. When solving ACOPF cases, \texttt{MIPS} has been shown to result in similar or fewer iterations than other solvers, such as \texttt{IPOPT} \cite{Baltzinger17}. In all of the cases tested in this paper, MIPS was able to converge starting both from the DC solution, a flat-start, and the learned solution. We also compare the reliability of our approach on the default interior point algorithm within the \texttt{fmincon} solver, which is more prone to convergence issues. For both solvers, the convergence criteria was set to the default of a $1e-6$ optimality gap, with a maximum iteration limit of 100.

\subsection{Hyperparameter tuning}
Using Scikit-learn's \texttt{GridSearchCV}, we varied the settings of the Random Forest and evaluated the performance of all combinations of settings to determine the best set of hyperparameters to include. We used $k = 3$ folds within the cross-validation, and the hyperparameters varied were the number of trees in the forest (\texttt{n\_estimators = [200, 300, 400, 500]}), the maximum number of levels in each tree (\texttt{max\_depth = [10, 15, 20]}), and the minimum number of samples required to split a node (\texttt{min\_samples\_split = [2, 3, 4, 5]}). The chosen hyperparameters are in Table I for each of the networks.

\begin{table}[t!]
\caption{Chosen hyperparameters for each of the networks.}
\centering
\begin{tabular}{|l|l|l|l|}
\hline
\textbf{Network} & \textbf{n\_estimators} & \textbf{max\_depth} & \textbf{min\_samples\_split} \\ \hline
\textbf{14-bus}  & 400                    & 15                  & 2                            \\ \hline
\textbf{57-bus}  & 200                    & 15                  & 4                            \\ \hline
\textbf{118-bus} & 400                    & 15                  & 2                            \\ \hline
\textbf{300-bus} & 400                    & 20                  & 3                            \\ \hline
\end{tabular} \label{tab:hyper}
\end{table}

\subsection{Prediction error}

The prediction error was normalized to represent a percentage and is defined as the absolute difference of the predicted variable $\hat{x}$ and the true optimum $x^*$ divided by the value of the true optimal variable value (i.e., the relative error):

\begin{align}
x_{err} = \frac{|\hat{x} - x^*|}{x^*}
\end{align}

\begin{figure}[t!]
    \includegraphics[width=0.5\textwidth]{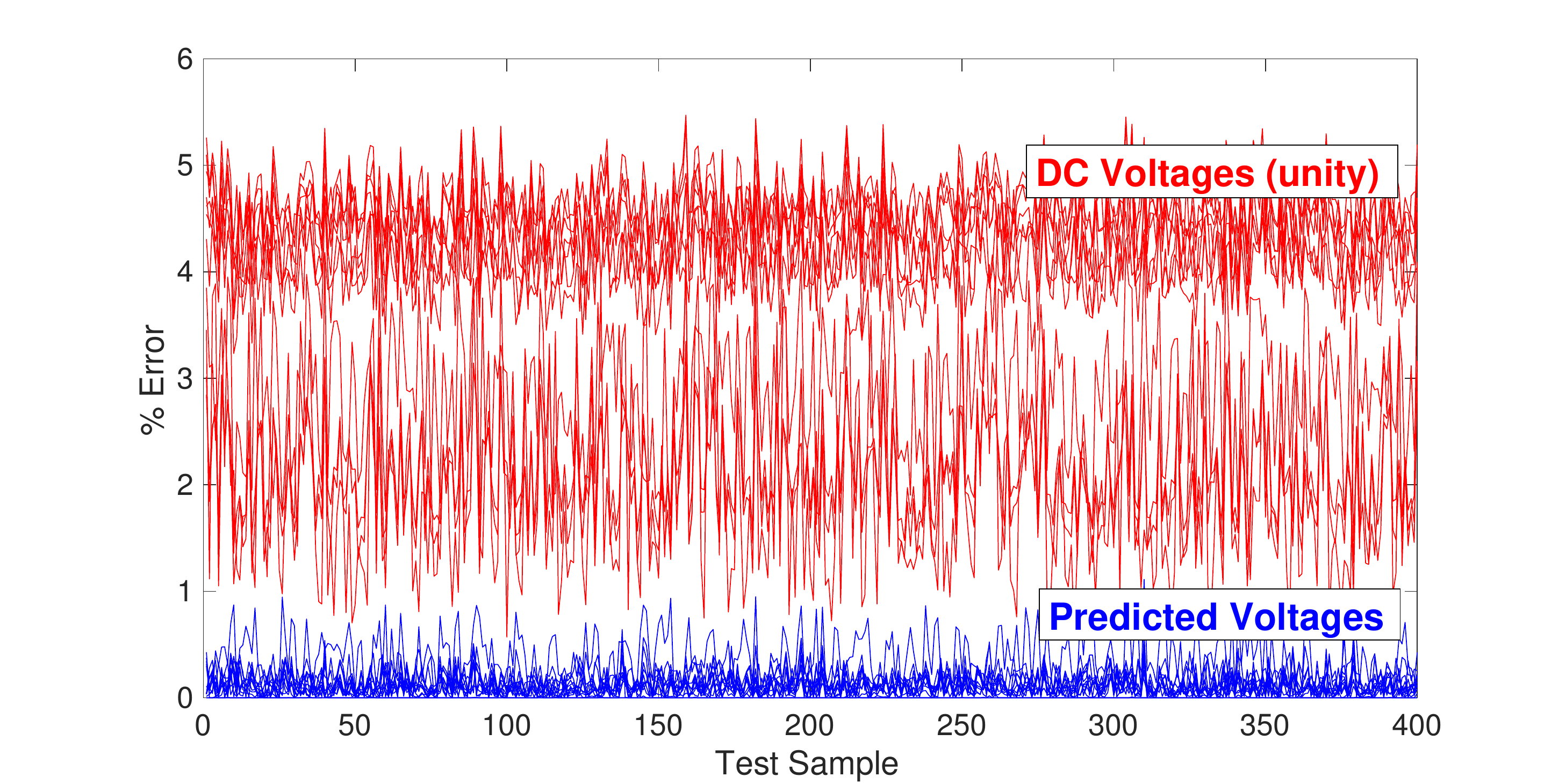}
    \caption{Relative error between predicted voltage magnitudes and optimal voltage magnitudes for each generator in the IEEE 14 bus test system. }
    \label{fig:per_voltage}
\end{figure}

\begin{figure}[t!]
    \includegraphics[width=0.5\textwidth]{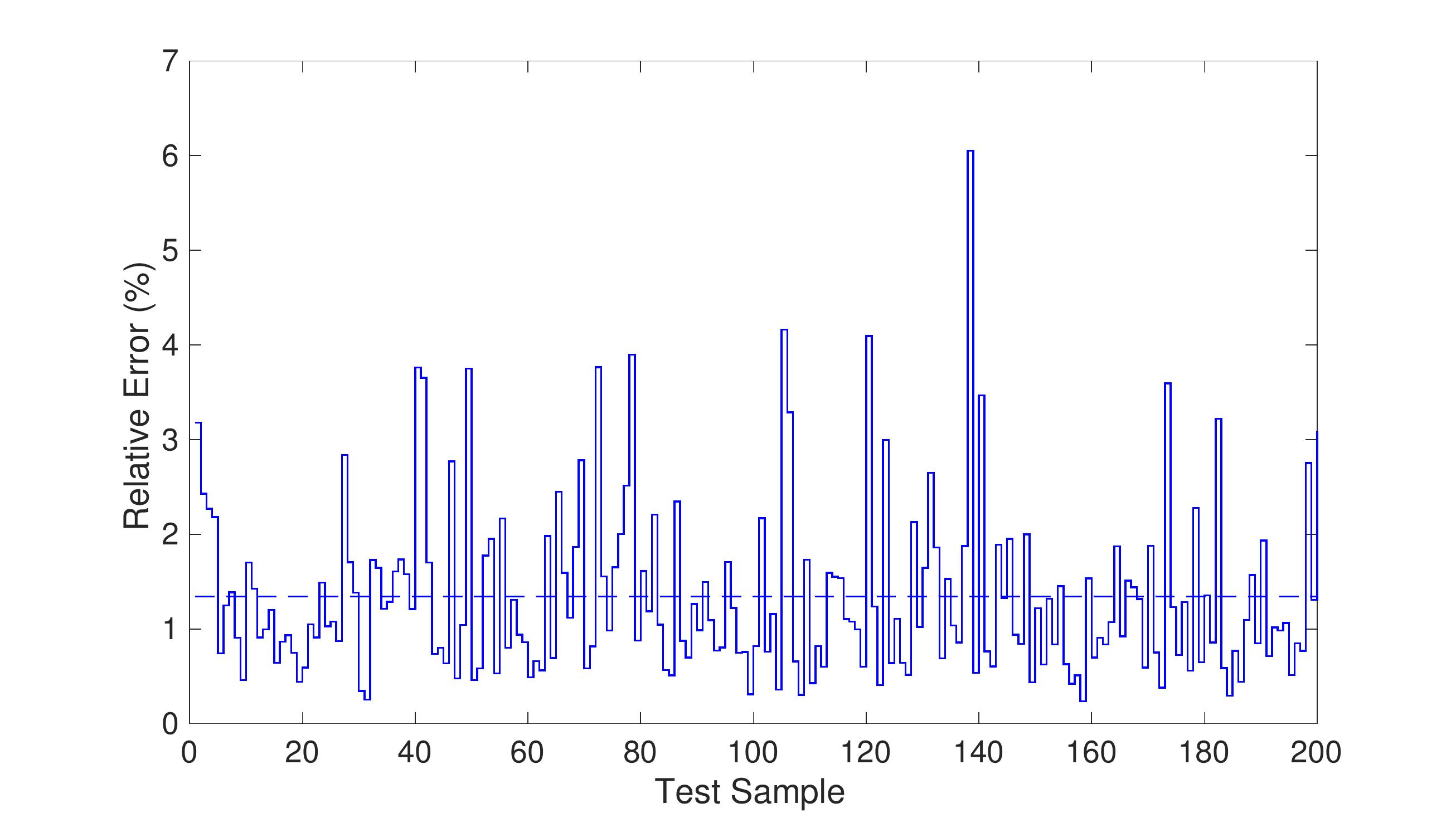}
    \caption{Relative error between predicted optimal generation and actual optimal generation for each bus in the IEEE 118 bus test system (truncated to 200 samples to show detail; dashed line is the mean of included points). }
    \label{fig:power_err}
\end{figure}

\noindent as shown for each bus voltage magnitude in the 14-bus network in Fig. \ref{fig:per_voltage}, the voltage magnitude predictions are very accurate, as this error is often less than 1\%. Table \ref{tab:error} shows the relative prediction errors for both voltage magnitude and optimal generation dispatch for each of the networks. for It is interesting to note that this is on par with or better than the error rates achieved by many AC power flow linearizations \cite{ SD_linear15, Bolognani15, Coffrin14}. From Fig. \ref{fig:per_voltage} we can see that the predicted voltages provide a much better initial guess than that provided by a DCOPF solution, which is often used as a traditional way of generating initial guesses. Figure \ref{fig:power_err} illustrates the accuracy of the optimal active power generation prediction on the 118-bus system by showing the average relative prediction error of all 54 generators in the network. The prediction is fairly accurate at estimating the optimal generation dispatch, and results in a very low error in most cases, and a lower error in all cases when compared to using the DC warm start. 

\begin{table}[]
\caption{Relative errors of predictions}
\centering
\begin{tabular}{|l|l|l|}
\hline
\textbf{Network} & \textbf{\begin{tabular}[c]{@{}l@{}}Avg. Relative\\ Error (Power)\end{tabular}} & \textbf{\begin{tabular}[c]{@{}l@{}}Avg. Relative \\ Error (Voltage)\end{tabular}} \\ \hline
\textbf{14-bus}  & 0.98\%                                                                         & 0.16\%                                                                            \\ \hline
\textbf{57-bus}  & 3.98\%                                                                         & 0.51\%                                                                            \\ \hline
\textbf{118-bus} & 2.12\%                                                                         & 0.01\%                                                                            \\ \hline
\textbf{300-bus} & 12.00\%                                                                        & 0.34\%                                                                            \\ \hline
\end{tabular} \label{tab:error}
\end{table}

\subsection{Improved computational time}

\begin{table}[t!]
\centering
\caption{Computational time required to predict an ACOPF solution given the network loads}
\begin{tabular}{|l|l|l|}
\hline
\textbf{Network} & \textbf{Number of Variables} & \textbf{Prediction Time} \\ \hline
\textbf{14-bus}  & 19                           & 9.04 ms                  \\ \hline
\textbf{57-bus}  & 64                           & 9.19 ms                  \\ \hline
\textbf{118-bus} & 172                          & 110.1 ms                 \\ \hline
\textbf{300-bus} & 369                          & 182.0 ms                 \\ \hline
\end{tabular} \label{tab:time_predic}
\end{table}

Training is performed offline, and obtaining a learned starting point was less than a second, meaning that the regression process does not significantly add to the overall computation time. In Table \ref{tab:time_predic}, the time required to predict the values with the Random Forest is shown. This is the time without distributed computations; i.e., without utilizing the benefits of parallelizing the Decision Tree regressions, meaning that this could potentially be performed even faster. Hence, approximate, yet accurate, solutions can be obtained in near real-time without even solving an optimization problem. 

\begin{table}[t!]
\caption{Total time to solve ACOPF on the testing dataset}
\centering
\begin{tabular}{|l|l|l|l|}
\hline
\textbf{Network} & \textbf{\begin{tabular}[c]{@{}l@{}}Tot. Time (s) \\ Learned\end{tabular}} & \textbf{\begin{tabular}[c]{@{}l@{}}Tot. Time (s)\\ DC\end{tabular}} & \textbf{\begin{tabular}[c]{@{}l@{}}Tot. Time (s)\\ Flat\end{tabular}} \\ \hline
\textbf{14-bus}  & 10.72                                                                     & 17.29                                                               & 14.06                                                                 \\ \hline
\textbf{57-bus}  & 15.18                                                                     & 24.76                                                               & 20.58                                                                 \\ \hline
\textbf{118-bus} & 25.46                                                                     & 33.49                                                               & 32.51                                                                 \\ \hline
\textbf{300-bus} & 56.68                                                                     & 67.82                                                               & 69.02                                                                 \\ \hline
\end{tabular} \label{tab:time}
\end{table}

When used to solve ACOPF on the testing dataset, we compared starting the \texttt{MIPS} solver with the learned warm start, a warm start from first solving the linearized DCOPF problem, and a flat start (voltage magnitudes set to 1.0 pu and angles and generation set to zero). From Table \ref{tab:time}, it is clear that the extra time required to solve a DCOPF in order to get an improved starting point over a flat start is generally not advantageous. Alternatively, using the Random Forest regression is fast enough that we see convergence time improvements over both the flat start and the DCOPF warm start.

\begin{figure}[t!]
    \includegraphics[width=0.5\textwidth]{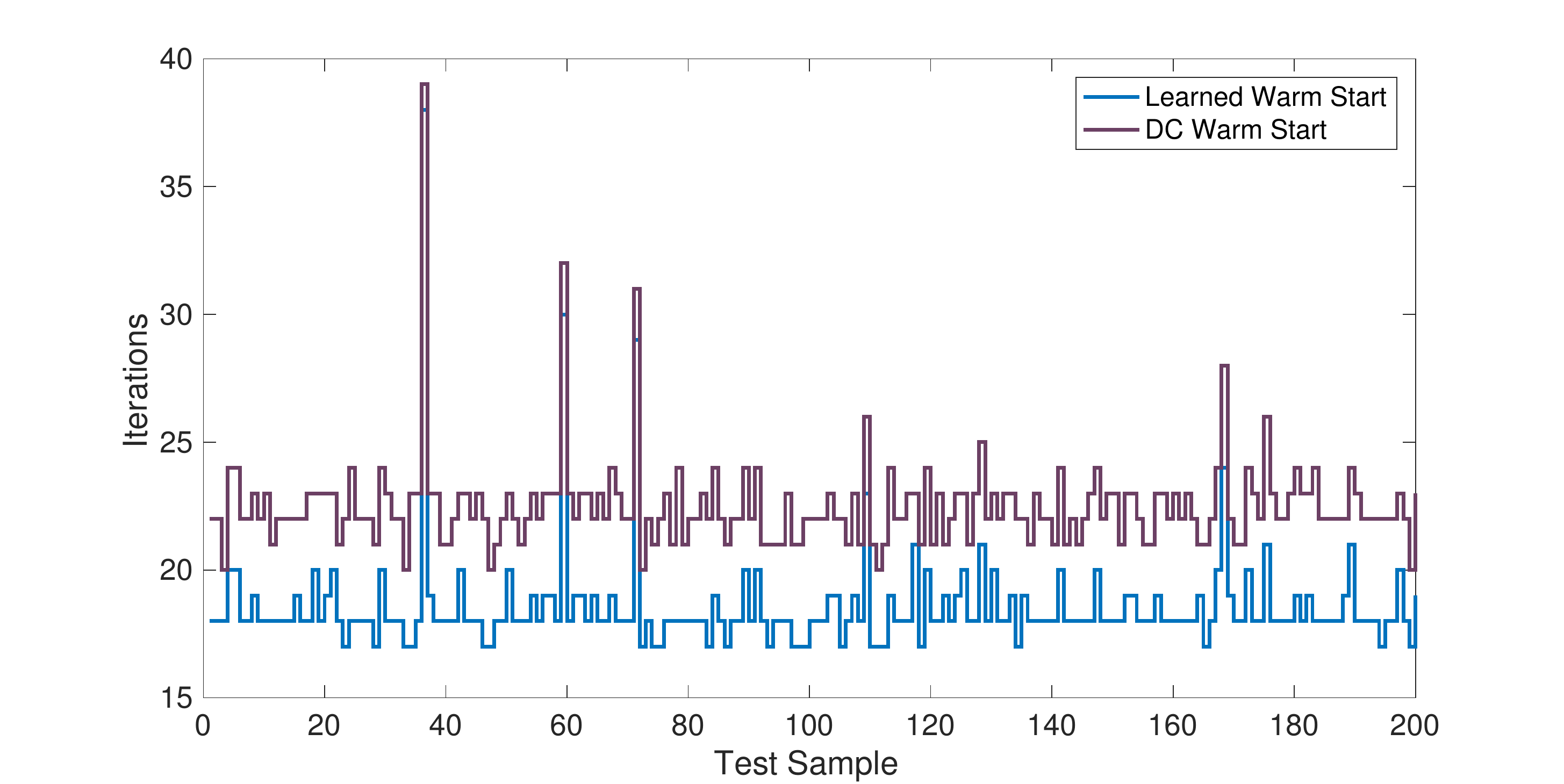}
    \caption{Number of iterations to convergence for solving ACOPF on the IEEE 300 bus network using both the learned warm start point and the DC optimal initial point. The learned initial point consistently outperforms a traditional warm start.}
    \label{fig:iter_300}
\end{figure}

As one example, in Fig. \ref{fig:iter_300}, the number of iterations to converge when using the learned start consistently outperforms the DCOPF solution in the largest test case. As stated in the previous section, in addition to lowering computational time, the learning approach can also provide a good proxy for finding approximate ACOPF solutions in real time, and its solution on par with or better than many linearization approaches. 

Due to the robustness of the \texttt{MIPS} solver, all of the cases, even with the DC and flat starting points, were able to converge and find a local minimum. However, we also decided to compare the ability of both the warm start and DC solution cases to successfully converge using less robust solvers (here, MATLAB's \texttt{fmincon}'s default interior-point solver). Using \texttt{fmincon}, convergence was not always achieved in either the learned start case or the DC start case. As Table \ref{tab:fmincon} illustrates, while the warm start provided a large benefit for successful convergence of the 57-bus case, it actually resulted in slightly worse performance for the 300-bus case. This is consistent with the findings in \cite{Venzke19}: the benefit of a warm start is solver-dependent. However, in the networks tested here, it is apparent that if a more robust solver is used, a learned warm start can provide benefits over a traditional warm start, and help convergence properties in difficult networks (i.e., the 57-bus).

\begin{table}[t!]
\caption{Percentage of runs that successfully converged using fmincon's interior point solver}
\centering
\begin{tabular}{|l|l|l|}
\hline
\textbf{Network} & \textbf{\begin{tabular}[c]{@{}l@{}}\% Converged\\ (Learned)\end{tabular}} & \textbf{\begin{tabular}[c]{@{}l@{}}\% Converged\\ (DC)\end{tabular}} \\ \hline
\textbf{14-bus}  & 100\%            & 100\%                                                                        \\ \hline
\textbf{57-bus}  & 100\%       & 53.00\%                                                                          \\ \hline
\textbf{118-bus} & 100\%           & 100\%                                                                      \\ \hline
\textbf{300-bus} & 99.50\%      & 100\%                                                                     \\ \hline
\end{tabular} \label{tab:fmincon}
\end{table}

\begin{figure}[t!]
    \includegraphics[width=0.5\textwidth]{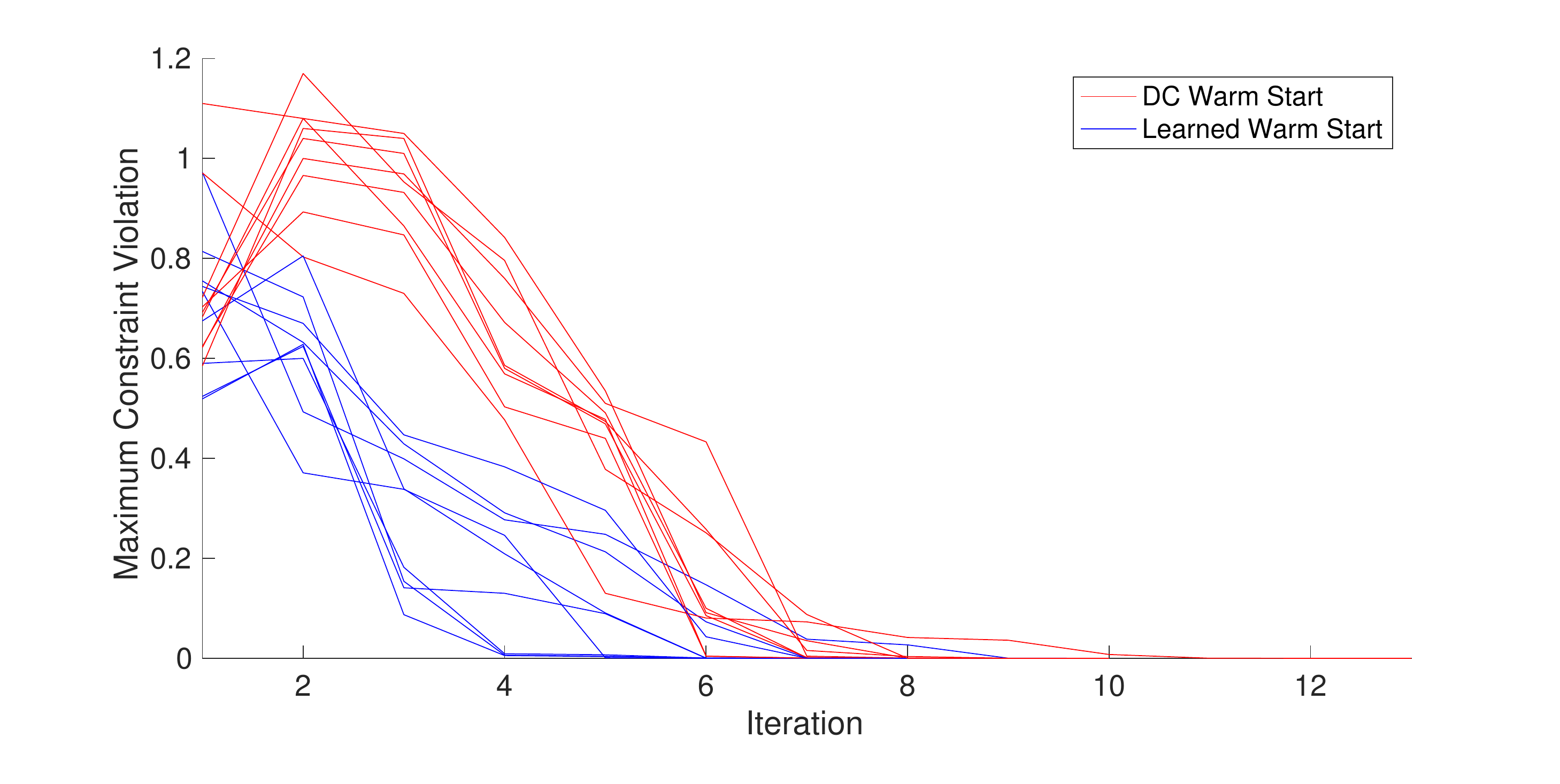}
    \caption{The measured maximum constraint violation for ten sample runs  on the IEEE 57-bus system using both a DC warm start and the learned initial point.}
    \label{fig:57_gap}
\end{figure}

In the presented results, almost every learned point resulted in an infeasible solution to the original ACOPF problem. However, as Fig. \ref{fig:57_gap} shows for 10 runs on the 57-bus system, while the solution is initially infeasible, the maximum constraint violation in the ACOPF is generally lower than that encountered when the DC initial point is used. 

\section{Conclusions and Future Work} \label{sec:conclusion}
We explored the benefits that using learning-based approaches can provide for determining a good starting point for solving ACOPF problems. In the considered networks, the learned warm start provided a faster convergence time over using a DC warm start or a flat start, and the overall benefits appear to be solver-dependent. We focused on Random Forests due to the distributed nature and thus computational benefits of such an approach. In addition, the accuracy of both predicted voltages and predicted optimal generation values is promising on its own for applications where finding approximate solutions quickly is of interest. Although generally infeasible, using the learned warm start solution generally results in lower maximum constraint violations than the DC initial point.

Future work will explore other machine learning techniques outside of Random Forests for predicting warm start solutions, such as neural networks or support vector regression. Iterations could perhaps be further lowered by including additional dependent variables in the prediction such as voltage angle and reactive power output; this and feature selection are important directions of future work. Methods to project infeasible solutions onto the feasible set or to ensure feasible predictions are of great interest. There was no sensitivity analysis performed on the sample size for the training set, which could also have an effect on the performance. 

\bibliographystyle{IEEEtran}
\bibliography{biblioNew.bib}



\end{document}